\documentclass[11pt]{amsart}

\usepackage{amscd}
\usepackage{amssymb}
\usepackage{boxedminipage}
\usepackage{enumerate}
\usepackage{euscript}

\theoremstyle{plain} \newtheorem{theorem}{Theorem}[section]
\theoremstyle{plain} \newtheorem{lemma}[theorem]{Lemma}
\theoremstyle{plain} \newtheorem{proposition}[theorem]{Proposition}

\newtheorem{corollary}[theorem]{Corollary}
\newtheorem{problem}{Problem}
\newtheorem{subproblem}[problem]{Subproblem}
\newtheorem{conjecture}[theorem]{Conjecture}
\newtheorem{conjecture*}[]{Conjecture}

\newcommand{\nr}{\refstepcounter{theorem}  
                   \noindent {\thetheorem .}}

\newcommand{\eks}{\medskip \noindent {\it Example \nr} }
\newcommand{\eksfin}{\medskip}
\newcommand{\rem}{\medskip \noindent {\it Remark \nr} }
\newcommand{\remfin}{\medskip}

\newcommand{\llabel}{\addtocounter{theorem}{-1}
\refstepcounter{theorem} \label}

\newcommand{\gL}{\mathcal{L}}

\newcommand{\gC}{\mathcal{C}}
\newcommand{\gO}{\mathcal{O}}

\newcommand{\EH}{{\mathbf H}}

\newcommand{\Hom}{\text{Hom}}

\newcommand{\Ext}{\text{Ext}}

\newcommand{\sus}{\subseteq}

\newcommand{\pil}{\rightarrow}

\newcommand{\mto}[1]{\stackrel{#1}\longrightarrow}

\newcommand{\impl}{\Rightarrow}

\newcommand{\iso}{\cong}
\newcommand{\te}{\otimes}

\newcommand{\tH}{\tilde{H}}

\newcommand{\la}{\lambda}

\newcommand{\dl}{\Delta}

\newcommand{\lk}{\text{lk}}
\newcommand{\lkd}{\lk_\Delta}

\newcommand{\del}{\Delta}
\newcommand{\delr}{\Delta_{-R}}
\newcommand{\dd}{{\dim \del}}

\renewcommand{\aa}{{\bf a}}
\newcommand{\bb}{{\bf b}}

\begin{document}

\title [Enriched homology and cohomology modules ]
{Enriched homology and cohomology modules of simiplicial complexes}
\author { Gunnar Fl{\o}ystad}
\address{ Matematisk Institutt\\
          Johs. Brunsgt. 12 \\
          5008 Bergen \\
          Norway}   
        
\email{ gunnar@mi.uib.no}

\begin{abstract}
For a simplicial complex $\del$ on $\{1,2, \ldots, n\}$ we define
enriched homology and cohomology modules. They are graded modules
over $k[x_1, \ldots, x_n]$ whose ranks are equal to the dimensions of
the reduced homology and cohomology groups.

We characterize Cohen-Macaulay, $l$-Cohen-Macaulay, Buchsbaum, and Gorenstein*
complexes $\del$,  
and also orientable homology manifolds in terms of the enriched modules. 
We introduce the notion of girth for simplicial complexes and make a conjecture
relating the girth to invariants of the simplicial complex.

We also put strong vanishing conditions on the enriched homology modules
and describe the simplicial complexes we then get. They are block
designs and include 
Steiner systems $S(c,d,n)$ and cyclic polytopes of even dimension.
\end{abstract}

\maketitle

{\Small MSC : 13F55, 55U10, 13D99.}

\section*{Introduction}
Given a simplicial complex $\dl$ on the set $[n] = \{1,2,\ldots, n\}$
and a field $k$, one has its reduced homology groups $\tilde{H}_i(\dl;k)$
which depend only on the topological realization
of $\dl$. However the combinatorial structure makes $\dl$ a richer
object than its topological realization. In this paper we define 
enriched homology modules $\EH_i(\dl;k)$ which are modules
over the polynomial ring $k[x_1,\ldots, x_n]$, denoted by $S$. 
They have the
property that the rank of $\EH_i(\dl;k)$ as an $S$-module is equal to
the dimension of $\tilde{H}_i(\dl;k)$ as a vector space over $k$, and
hence may be considered as an enrichment of the reduced homology 
groups.

The enriched homology modules are defined as the homology modules of
the cellular complex, see \cite{BS}, associated to $\dl$ by attaching
the variable $x_i$ to the vertex $i$. Similarly we may define the
enriched cohomology modules $\EH^i(\dl;k)$ as the cohomology modules
of the cocellular complex. 

The classical criterion of Hochster for when a simplicial complex is
Cohen-Macaulay, via the reduced homology of its links, translates in this
setting to the criterion that the top enriched cohomology module
$\EH^{\dim \dl}(\dl;k)$ is the only non-vanishing cohomology module.
We also give a criterion for when $\dl$ is Buchsbaum via its 
enriched cohomology modules.

The enriched cohomology modules turn out to be Alexander duals of 
the Matlis duals of the local cohomology modules of the Stanley-Reisner
ring. Hence they contain exactly the same information as these local
cohomology modules. The enriched cohomology modules present
both another conceptual approach to this information, and the
information appears in quite distinct algebraic forms. Therefore algebraic
questions which are not interesting for local cohomology modules
or their Matlis duals, turn out to be interesting for enriched cohomology
modules.

\medskip
For instance when $\dl$ is Cohen-Macaulay we 
show in \cite{Fl} that the top cohomology module
can occur as an $l-1$'th syzygy module in an $S$-free resolution if
and only if $\dl$ is $l$-Cohen-Macaulay as defined by Baclawski \cite{Ba}.
(This notion means that $\dl$ restricted to each subset $R$ of $[n]$ of 
cardinality $n+1-l$ is a Cohen-Macaulay complex of the same dimension as
$\dl$.) 

For $l$-Cohen-Macaulay simplicial complexes, the dimension of $\del$
is less or equal to $n-1-l$ (unless $\del$ is the $n-l$-skeleton of
the simplex), and the girth of $\del$, a suitable generalization
of the well-known notion for graphs, is less or equal to $n+2-l$.
We investigate $l$-Cohen-Macaulay simplicial complexes with one or 
more of these
extremal values. In particular we make a conjecture concerning invariants of
$l$-Cohen-Macaulay complexes of maximal girth, and prove the conjecture
in the case where the dimension is submaximal equal to $n-1-l$.

Next we consider the case when the top cohomology module may be identified
as an ideal in the polynomial ring $S$. 
When the top cohomology module is the only non-vanishing one, i.e. when $\del$
is CM, we show that this happens exactly 
when $\dl$ is a Gorenstein* simplicial complex. More generally for 
Buchsbaum complexes of positive dimension, 
we show that the top cohomology module is an ideal
if and only if the complex is a connected orientable homology manifold.

\medskip
In the end we put strong vanishing conditions on the enriched
homology and cohomology
modules and investigate what kind of simplicial complexes this corresponds to. 
For an $l$-Cohen-Macaulay $\dl$ 
there is only one non-vanishing cohomology module.
We require that there also be at most one
non-vanishing homology module $\EH_i(\dl;k)$ for $i < \dim \dl$
(when $l \geq 2$, $\EH_{\dim \dl}(\dl;k)$ does not vanish), and that
$\dl$ has maximal girth, which is $n+2-l$. 
In \cite{FV} we introduced the notion of $\dl$ being bi-Cohen-Macaulay,
meaning that both $\dl$ and its Alexander dual are Cohen-Macaulay.
We show that the above vanishing condition is equivalent
to the condition that the restriction of $\dl$ to each subset $R$ of $[n]$ 
of cardinality $n+1-l$ is bi-Cohen-Macaulay of the same dimension and 
frame dimension
as $\dl$. (The frame dimension is by definition the dimension of the
maximal complete skeleton of the simplex on $[n]$, which is contained in 
$\dl$.) In this case the $f$-vector of $\dl$ is completely determined
by $l$, $n$, the dimension, and frame dimension of $\dl$, and
we call $\dl$ an {\it $l$-Cohen-Macaulay design}.

As examples, $1$-Cohen-Macaulay designs are exactly the bi-Cohen-Macaulay
simplicial complexes. $l$-Cohen-Macaulay designs of submaximal dimension
$n-l-1$ are exactly the Alexander duals of Steiner systems $S(l-1,m,n)$
(where $m$ is related to the frame dimension of $\dl$). $2$-Cohen-Macaulay
designs of dimension twice the frame dimension pluss one, are examplified
by the boundaries of cyclic polytopes of even dimension.

We also introduce the more general class of $(a,b)$-designs which 
specialize to $a\!+\!1$-CM designs when $b=0$, and show that they
are block designs. As examples 
$2$-neighbourly triangulations of surfaces are $(1,1)$-designs, 
and more generally
homology manifolds with certain extremal behaviour of the Euler 
characteristic are $(1,1)$-designs.

\medskip

The organization of the paper is as follows. In Section $1$ we define
the enriched homology and cohomology modules, and recall some basic facts
on cellular complexes, Koszul duality and square free modules (as defined
by Yanagawa \cite {Ya1}). We translate the link criterion of Hochster for
$\dl$ being Cohen-Macaulay, and the criterion of Schenzel for $\dl$ being
Buchbaum, to the setting of cohomology modules. We also translate the
restriction criterion of Hochster for $\dl$ being Cohen-Macaulay
to the setting of homology modules.

In Section $2$ we 
first recall a result from \cite{Fl} providing motivation for the 
usefulness of enriched cohomology. Namely that the top cohomology
module of CM complexes can occur as an $l-1$'th syzygy module
in an $S$-free resolution iff $\dl$ is $l$-Cohen-Macaulay.
Then we go on to investigate $l$-CM simplicial complexes of maximal
girth. In particular we make a conjecture concerning the invariants
of such and prove this when the complex has submaximal dimension.

In Section $3$ we consider when the top cohomology module is
an ideal in the polynomial 
ring. For Cohen-Macaulay complexes 
this happens when it is Gorenstein*. More generally,
for Buchsbaum complexes this happens when it is 
an orientable homology manifold.

In Section $4$ we put strong vanishing conditions on the homology
and cohomology modules and investigate the notion of $l$-Cohen-Macaulay 
designs. We also introduce the more general notion of $(a,b)$-designs.
This section also contains many examples of such designs.
Finally in Section $5$ we give some problems and conjectures.

\medskip
\noindent {\bf Note.} This paper is to a large extent a complete rewriting
of a previous preprint, { ``Hierarchies of simplicial complexes via the
BGG-correspondence''}. Also Proposition \ref{HomProCoh} and \ref{TopTheGor} 
have been generalized to
cell complexes in \cite{Fl}.

\medskip
\noindent {\bf Acknowledgements.} I thank an anonymous referee for 
suggestions simplifying the proofs of Theorems \ref{TopTheGor}
and \ref{TopTheBb}, and in general for improving the exposition
of the paper.

\section{Enriched homology and cohomology modules}

\subsection{Notation} Denote by $[n]$ the set $\{1,2, \ldots, n \}$. 
A simplicial complex
$\del$ on $[n]$ is a family of subsets of $[n]$ such that if $F$ is 
in $\del$ and $G$ is a subset of $F$, then $G$ is in $\del$.

An element $F$ in $\del$ is called a face of $\del$. If $F$ has cardinality
$f$, its dimension is $f-1$. If $d$ is the maximal cardinality of a face,
the {\it dimension} of $\del$ is $d-1$. A face of maximal cardinality 
is called a facet.
The maximum $i$ such that all $i$-sets are in $\del$ is denoted
by $c$ and we call $c-1$ the {\it frame dimension} of $\del$.

If $R$ is a subset of $[n]$ denote by $\del_R$ the restricted simplicial
complex on $R$, 
i.e. consisting of the faces $F$ contained in $R$. It will also
be convenient to have the notation $\delr$ for the restriction to 
the complement $[n]\backslash R$. 

For $Q$ a subset of $[n]$ the link of $Q$, $\lkd Q$, is the simplicial complex
on $[n]\backslash Q$ consisting of the subsets $F$ of $[n]\backslash Q$
such that $Q \cup F$ is a face of $\del$.

A convention we will often use is that if a set is denoted by an
upper case letter, say $R$, then the lower case letter $r$ will denote
the cardinality of the set $R$.

\subsection{Enriched homology modules}
Given a field $k$, one has the {\it augmented oriented chain complex}
$\tilde{\gC}(\del;k)$. The group $\tilde{\gC}_i(\del;k)$ is the vector space
$\oplus kF$ with basis consisting of the faces of dimension $i$, and 
the differential defined by
\[ F \mapsto \sum_{\dim F^\prime = i-1} \epsilon(F,F^\prime) F^\prime. \]
where $\del \times \del \mto{\epsilon} \{-1,0,1\}$ is a suitable
incidence function (see \cite{BH}). The homology groups 
$\tilde{H}_i(\del;k)$ of this complex
are the {\it reduced homology groups} of $\del$ and depend only on the
topological realization of $\del$. 

However the combinatorial structure makes the simplicial complex a richer
structure than its topological realization. We define enriched homology 
modules of the simplicial complex as follows. Let $S$ be the 
polynomial ring $k[x_1, \ldots, x_n]$. 
We get a complex $\gL(\del;k)$ of free $S$-modules
by letting $\gL_i(\del;k)$ be the free $S$-module $\oplus SF$ with
generators the faces of dimension $i$
and differential
\[ F \mapsto \sum_{F = F^\prime \cup \{l\}} \epsilon(F,F^\prime)
x_l \cdot F^\prime. \]
We define the {\it enriched homology modules} $\EH_i(\del;k)$ 
(or just $\EH_i(\del)$) to 
be the homology modules of this complex.

 There are two sources of inspiration for this definition. First there is
the theory of cellular complexes developed by Bayer, Peeva, and Sturmfels,
\cite{BS}, \cite{BPS}. The complex $\gL(\del;k)$ is the cellular complex
obtained by attaching the monomial consisting of the single variable
$x_i$ to vertex $i$. Another approach comes from the Koszul duality between
the symmetric algebra $S$ and the exterior algebra $E$ in $n$ variables,
and how this relates the module categories of these algebras, (see
\cite {EFS} and \cite{FV} for some recent articles). We explain this in 
some detail. This will also make the incidence function completely explicit.

\medskip

Let $V$ be the vector space on generators $e_1, \ldots, e_n$ and 
\[ E = E(V) = \oplus_{i = 0}^{n} \wedge^i V\]
be the exterior algebra.
We let $W = V^*$ be the dual vector space with dual basis $x_1, \ldots, x_n$
and identify the polynomial ring $S$ as the symmetric algebra Sym$(W)$.
We consider $V$ to have degree $-1$ and $W$ to have degree $1$.
For a graded module $M$ over $E$ or $S$ we denote by $M(j)$ the module
shifted $j$ steps to the left i.e. $M(j)_d = M_{j+d}$.

Any graded (left) $E$-module $M$ gives rise to a complex of $S$-modules
\begin{equation} \label{HomLabLsek}
L(M) : \, \,  \cdots \pil S \te_k M_i \mto{d_i} S \te_k M_{i-1} \pil
\cdots 
\end{equation}
where $d_i$ sends $s \te m$ to $\sum_j sx_j \te e_jm$. Note that the 
degree of $s \te m$ is the sum of the degrees of $s$ and $m$.

Given the simplicial complex $\del$, we can form the {\it exterior face ring}
$k\{\del\}$, which is the quotient of $E$ by the monomial ideal $J_\del$ 
consisting of monomials $e_{i_1} \cdots e_{i_r}$ such that 
$\{i_1, \ldots, i_r\}$ is not a face of $\del$.
Let $C_\del$ be the graded dual vector space 
of the exterior face ring $k\{\del\}$. It is a module over the exterior face 
ring. As a vector space it has a basis consisting of monomials 
$x_{i_1} \cdots x_{i_r}$ where
$\{i_1, \ldots, i_r\}$ range over the faces of $\del$. Left multiplication
with $e_1 + e_2 + \cdots +e_n$ gives a differential $d$ on $C_\del$ and
the reduced homology of $\del$ is given by
\[ \tilde{H}_p(\del,k) = H_{p+1}(C_\del,d). \] 
The enriched homology of $\del$ is given by 
\[ \EH_p(\del;k) = H_{p+1} (L(C_\del)) .\]
We denote $L(C_\del)$ by $L(\del;k)$ (or just $L(\del)$). 
Note that compared to $\gL(\del;k)$ it is shifted one step to the left.

\medskip

If ${\bf b}$ is a multidegree in 
${\bf N}^n$, the 
{\it support} of ${\bf b}$ is the set of non-vanishing coordinates.
The following explicitly describes the multigraded parts of the homology 
module.

\begin{lemma} \label{HomLemMul} 
For a multidegree ${\bb}$ in ${\bf N}^n$, let $R$ be its support.
Then
\[\EH_p(\del)_\bb \iso \tH_p(\del_R). \] 
\end{lemma}

\begin{proof} This follows from the above description of $L(\del)$ as
a cellular complex and the ideas in the proof of Proposition 1.1 in \cite{BS}.
See also Corollary \ref{HomCorResR} below.
\end{proof}

In particular we shall consider $\tH_p(\del_R)$ to have multidegree $R$.

\rem From the above lemma and Hochster's description of the resolution of the
Stanley-Reisner ring $k[\del]$, see \cite{St} or originally \cite{Ho},
we see that the homology module $\EH_p(\del)$ corresponds to the
$p+1$'th linear strand in the resolution. Thus the collection of homology
modules is equivalent to give the linear strands of the resolution of 
the Stanley-Reisner ring. Our approach  gives another point of view
to this and new questions are natural to ask.

\subsection{Square free modules}
The notion of square free $S$-module was introduced by Yanagawa \cite{Ya1}.
A ${\bf N}^n$-graded $S$-module $M$ is {\it square free} if 
$M_\bb \mto{x_i} M_{\bb + u_i}$ is an isomorphism for every 
$\bb$ in ${\bf N}^n$ with $i$ in the support of $\bb$, and
where $u_i$ is the $i$'th coordinate vector.

It follows from the description of $L(\del)$ as a cellular complex, \cite{BS},
that the enriched homology modules
are square free $S$-modules. Part a. and b. in the
following proposition are quotes and consequences of \cite{Ya1}, 
Lemma 2.2, Corollary 2.4, and Proposition 2.5. 

\begin{proposition} \label{HomProSQ} Let $M$ be a square free $S$-module.

a. The minimal prime ideals are the $(x_i)_{i \not \in R}$ where
$R$ are the maximal subsets of $[n]$ with $M_R$ nonzero.

b. For $R$ maximal as above let $P$ be the corresponding prime ideal. There
is a natural map ($M_R$ has multidegree $R$)
\[ S/P \te_k M_R \pil M \]
which becomes an isomorphism when localized at $P$.

\end{proposition}




The following is the justification for the terminology enriched homology
modules.

\begin{corollary} The $S$-module rank of $\EH_p(\del)$ is the dimension of 
$\tH_p(\del)$ as a vector space over $k$.
\end{corollary} 

\begin{proof} Let $M$ be $\EH_p(\del)$ and $R = [n]$. 
Proposition \ref{HomProSQ} b., gives a natural map
$S \te_k \tH_p(\del) \pil \EH_p(\del)$ which becomes an isomorphism when
localized at $(0)$.
\end{proof}

There is also a notion of square free modules over the exterior algebra
$E = E(e_1, \ldots, e_n)$, namely a multigraded module $M$ is called  
square free if $M_{\bb}$ is nonzero only if $\bb$ is a characteristic
vector of some $R \sus [n]$, i.e. a vector such that $b_i$ is $1$ for
$i$ in $R$ and $b_i$ is zero otherwise.
(This notion is a variation of the one defined by T.R\"omer, \cite{Ro}.
This is because we consider the $e_i$ to have negative degrees. 
According to our convention
$E$ is not square free, but its dual $E(W)$ is.)
We may note that there is an equivalence of categories between the 
square free modules over $S$ and $E$, \cite{Ro}.

For a square free module $M$ over $E$, denote by $M_{|R}$ the restriction
to $E_R = E(e_i, i \in R)$, i.e. $(M_{|R})_{\bb}$ is $M_{\bb}$ 
when the support of $\bb$ is in $R$, and zero otherwise.
For a square free module over $S$ we may restrict it to
$S_R = k[x_i, i \in R]$, the definition is by the same formulation as above.
Now as in (\ref{HomLabLsek}), define the functor $L_R$ on a module $M^\prime$ over $E_R$ to be
\[ L_R(M^\prime) : \, \,  \cdots \pil S_R \te_k M^\prime_i 
\mto{d_i} S_R \te_k M^\prime_{i-1} \pil
\cdots  .\]

\begin{lemma} \label{HomLemRes} Let $M$ be a square free module over $E$.

a. $L(M)_{|R} = L_R(M_{|R})$.

b. $H^p(L(M))_{|R} = H^p(L_R(M_{|R}).$
\end{lemma}

\begin{proof}
a. This is straightforward. b. Follows since restriction is an exact functor.
\end{proof}

\begin{corollary} \label{HomCorResR} 
$\EH_p(\del)_{|R} = \EH_p(\del_R)$.
\end{corollary}

\begin{proof}
This is because $(C_\del)_{|R} = C_{\del_R}$.
\end{proof}

\subsection{Alexander duals}

The Alexander dual $\del^*$ of the simplicial complex $\del$, is the simplicial
complex on $[n]$ consisting of the subsets $F$ such that $F^c$ is not in 
$\del$.

The Alexander dual of the square free module $M$ over $E(e_1, \ldots, e_n)$ is
the square free module $M^\vee = \Hom_k(M,k(-[n]))$.
Note that these notions are related by an exact sequence
\[ 0 \pil C_{\del^*} \pil E(W) \pil (C_\del)^\vee \pil 0. \]

Let $M$ be a square free module over $S$. Via the
equivalence of categories this corresponds to a square free module
over $E$. We may take the Alexander dual of this module, and via the equivalence
we again get a square free module $M^*$ over $S$, called the Alexander dual
of $M$, \cite{Mi} and \cite{Ro}. Explicitly,
for a subset $F$ of $[n]$, $(M^*)_F$ is $\Hom_k(M_{F^c},k)$ where
$F^c$ is the complement of $F$ in $[n]$. If
$i$ is not in $F$, the multiplication map
\[ (M^*)_F \mto{x_i} (M^*)_{F \cup \{i\}} \]
is the dual of the multiplication map
\[M_{F^c \backslash \{i \}} \mto{x_i} M_{F^c}. \]
By obvious extension this defines $(M^*)_\bb$ for all $\bb$ in ${\bf N}^n$ and
all multiplications by variables.

\subsection{Enriched cohomology modules}

The reduced cohomology groups $\tH^p(\del;k)$ of the simplicial complex
$\del$ are the cohomology groups of the dualized complex
\[\Hom_k(\tilde{\gC}(\del;k), k). \]

We define the {\it enriched cohomology modules} $\EH^p(\del;k)$
(or just $\EH^p(\del)$) as the cohomology modules of the 
dualized complex 
\[\gL(\del;k)^\vee = \Hom_S(\gL(\del;k), \omega_S) \]
where $\omega_S$ is the canonical module of $S$,
isomorphic to $S(-{\bf 1})$ where ${\bf 1}$ is the multidegree
$(1,1,\ldots,1)$.

In terms of the Koszul duality correspondence this is given as follows.
The reduced simplicial cohomology group is
\begin{equation} 
 \tH^p(\del;k) = H^{p+1}((C_\del)^*). \label{HomForSiKo}
\end{equation}

The complex $\gL(\del;k)^\vee$ identifies as the complex 
$L((C_\del)^\vee)[-n]$ 
(here $[-n]$ denotes the complex shifted $n$ steps to the right)
except that the former is shifted one step to the left. The 
enriched  cohomology modules are
\[ \EH^p(\del;k) = H^{p+1}(L((C_\del)^\vee)[-n]).\]

The following describes these cohomology modules in greater
detail.

\begin{proposition} \label{HomProCoh}
$\EH^p(\del)_{|R} = \EH^{p-r^c}(\lk_{\del} R^c)$. In particular 
$\EH^p(\del)_{\bb}$ is isomorphic to $\tH^{p-q}(\lkd Q)$
where $Q$ is the complement in $[n]$ of the support of $\bb$.
\end{proposition}

\begin{proof}
This follows by Lemma \ref{HomLemRes} since
$(C_\del)^\vee_{|R}$ equals $(C_{\lk_\del R^c})^\vee(-R^c)$.
\end{proof}

The above corollary suggests that enriched cohomology modules
are related to the local cohomology modules of the Stanley-Reisner ring.
In fact, they contain exactly the same information as the following
shows.

 \begin{theorem} \label{HomTheExt}
$\EH^p(\del)$ is the Alexander dual of $\Ext_S^{n-p-1}(k[\del], \omega_S)$,
which again is Matlis dual to the local cohomology module
$H^{p+1}_m(k[\del])$.
\end{theorem}

\begin{proof}
The last statement is local duality. The former follows by \cite[Thm. 5.6.3]{BH}
or alternatively from \cite[Prop. 3.1]{Ya1}.
\end{proof}


\subsection{Cohen-Macaulay and Buchsbaum simplicial complexes}

A simplicial complex $\del$ is called Cohen-Macaulay (CM) 
if the Stanley-Reisner
ring $k[\del]$ is a Cohen-Macaulay ring. By a criterion of Hochster, \cite{Ho}
or \cite[II.4]{St}, this is equivalent to
\begin{equation} \label{HomForHoL} \tH_p(\lkd R) = 0 \mbox{ for }
p+r < \dim \del. \end{equation}
This gives the following criterion.

\begin{proposition} \label{HomProCM}
$\del$ is Cohen-Macaulay iff the cohomology modules $\EH^p(\del)$
vanish for $p < \dim \del$, or alternatively $L(\del)^\vee$
is a resolution of $\EH^{\dim \del}(\del)$.
\end{proposition}

\begin{proof}
Since $k$ is a field, $\tH^p(\lkd R)$ is isomorphic to $\tH_p(\lkd R)^*$.
By Corollary \ref{HomProCoh} and Hochster's criterion (\ref{HomForHoL})
this translates exactly to the above statement.
\end{proof}

This shows that the well-known link criterion for $\del$ being 
Cohen-Macaulay
is encoded quite compactly in the cohomology modules of $\del$. 

\rem In general, when $k$ is not a field, 
the condition of only one non-vanishing cohomology module is weaker
then being Cohen-Macaulay, only implication to the right holds
in Proposition \ref{HomProCM}.

\begin{proposition} When $\dl$ is CM, the top cohomology module
is the Alexander dual of the canonical module $\omega_{k[\dl]}$ of
the Stanley-Reisner ring.
\end{proposition}

\begin{proof} By Theorem \ref{HomTheExt}, $\EH^{\dim \dl}(\dl)$ is 
Alexander dual to 
\[\Ext_S^{n-\dim \dl -1}(k[\dl],\omega_S)\] 
which is the canonical module of the Stanley-Reisner ring.
\end{proof}

The simplicial complex $\del$ is called Buchsbaum if $k[\del]$ is a 
Buchsbaum ring. By a criterion of Schenzel, \cite{Sc} or \cite[II.8]{St},
this is equivalent to 
\begin{equation} \label{HomForBub} \tH_i(\lkd R) = 0
\mbox{ for } i+r < \dim \del, \, r \geq 1. \end{equation}
This gives the following criterion.

\begin{proposition} \label{HomProBub}
$\del$ is Buchsbaum iff 
\[ \EH^p(\del) \iso S \te_k \tH^p(\del) \mbox{ for } p < \dim \del. \]
(Note that $\tH^p(\del)$ has multidegree ${\bf 1}$, see convention
after Lemma \ref{HomLemMul}.)
\end{proposition}

\begin{proof} When $p < \dim \del$
the criterion (\ref{HomForBub}) says that $\EH^p(\del)_\bb$ vanishes unless
the support of $\bb$ is the whole of $[n]$. But any square free module $M$
with this property must be of the form $S \te_k M_{[n]}$.
\end{proof}

An alternative criterion for the  Cohen-Macaulayness of $\del$ is the
following, \cite[p.197]{Ho},
\begin{equation} \label{HomForHoR}\tH_p(\del_{-R}) = 0 \mbox{ for }
 p+r < \dim \del. \end{equation}
This follows directly from Hochster's formula for the
multigraded Tor's in the resolution of the Stanley-Reisner ring.

For a square free module $M$, the 
codimension of $M$ is $\geq c$ iff
$M_{[n]\backslash F}$ is zero for all
$F$ of cardinality less than $c$. 
Hochster's criterion (\ref{HomForHoR}) is then equivalent to the following.

\begin{proposition} \label{HomProHo}
$\del$ is CM iff each homology module
$\EH_{\dim \del - i}(\del)$ has codimension $\geq i$.
\end{proposition}

\section{$l$-Cohen-Macaulay simplicial complexes}
\label{TopSek}

For Cohen-Macaulay simplicial complexes there is only one nonvanishing
cohomology module, $\EH^{\dim \del}(\del)$. It is therefore natural to put
various algebraic conditions on this module and investigate what
properties this corresponds to for the simplicial complex.

\subsection{$l$-CM simplicial complexes}
In \cite{Ba}, K. Baclawski introduced the notion of {\it $l$-Cohen-Macaulay} 
simplicial complexes, which geometrically corresponds to higher connectivity. 
The simplicial complex $\del$ is said to be $l$-Cohen-Macaulay ($l$-CM)
if $\del_{-R}$ is Cohen-Macaulay of the same dimension as $\del$, 
for all $R$ of cardinality $\leq l-1$.
For instance if $\del$ is a graph, then $\del$ is $l$-CM iff it is 
(vertex) $l$-connected.

In \cite{Fl} we prove that this property has nice descriptions in terms
of the top cohomology module and also in terms of the homology modules,
generalizing Propositions \ref{HomProCM} and \ref{HomProHo}. 

\begin{theorem}[\cite{Fl}] \label{TopTheSyz}
The following are equivalent for a simplicial complex $\del$.

\noindent 1. $\del$ is CM and $\EH^{\dd}(\del)$ 
can occur as an $l-1$'th syzygy 
module in an $S$-free resolution.

\noindent 2. The codimension of $\EH_{\dd-i}(\del)$ is greater or equal to 
$(l-1)+i$ for $i \geq 1$.

\noindent 3. $\del$ is $l$-CM.
\end{theorem}

\rem In the theory of polytopes, Balinski's theorem, \cite{Zi}, says that the
$1$-skeleton of a $d$-dimensional polytope is $d$-connected. 
In \cite{Fl} we show a comprehensive generalization of 
Balinski's theorem, namely that the codimension $r$ skeleton of an $l$-CM 
simplicial complex is $l+r$-CM. This is a rather immediate consequence of the
above theorem.

\rem From \cite{Ba} we mention the following two properties
which produces new $l$-CM complexes.
i) If $\del$ is an $l$-CM simplicial complex, then any link $\lkd Q$ is 
$l$-CM.
ii) If $\del_1$ and $\del_2$ are $l$-CM simplicial complexes,
their join $\del_1 * \del_2$ is $l$-CM. In particular since 
$l$ vertices is $l$-CM, if $\del$ is $l$-CM, 
the $l$-point suspension of $\del$ is $l$-CM.

\subsection{Maximal girth and maximal dimension}

For a simplicial complex $\del$ we define its {\it girth} to be 
the smallest degree in which the top homology module $\EH_{\dd}(\del)$ is 
nonzero. Since homology modules are square free the girth is $\leq n$
provided the top homology module does not vanish. If it vanishes we define
the girth to be $n+1$.

If $\del$ is a graph, this specializes to the notion of girth for graphs,
the length of a cycle of minimal length.

\begin{proposition} \label{TopProGiSu}
Let a non-empty $\del$ be $l$-CM 

a. Its girth $\leq n+2-l$.

b. The cardinality $d$ of a facet is $\leq n-l$, unless $\del$ is the 
$n-l$-skeleton of the $n-1$-simplex.
\end{proposition}

\begin{proof}
a. We want to show that $H_{\dd}(\del_{-R})$ is non-zero for some $R$
of cardinality $l-2$. By restricting to $\del^\prime = \del_{-R}$ for some
$R$ of cardinality $l-2$, it will be sufficient to show this for $2$-CM 
$\del^\prime$. 
Now it is easy to see that the top cohomology module of any simplicial
complex is non-zero unless the simplicial complex is empty. 
For instance this may be seen from Proposition \ref{HomProCoh}
by taking $Q$ to be a facet of the simplicial complex.
But then applying Theorem \ref{TopTheSyz} 
to the $2$-CM $\del^\prime$, its non-zero
enriched top cohomology module $\EH^{\dim \del^\prime}(\del^\prime)$
must have rank $\geq 1$ as an $S$-module. Then $\tH^{\dim \del^\prime}
(\del^\prime)$ is nonzero and so also the reduced homology module
$\tH_{\dim \del^\prime}
(\del^\prime)$.

b. The restriction $\del_{-R}$ has the same dimension as $\del$ for all $R$
of cardinality $l-1$. Hence $d \leq n+1-l$. If $d$ is $n+1-l$, each 
$\del_{-R}$ would be a simplex, since $[n]\backslash R$ has cardinality 
$n+1-l$.
But then $\del$ would be the $n-l$-skeleton of the simplex on $[n]$.
\end{proof}

\rem If $\dl$ is $l$-CM, its codimension one skeleton will have girth
$\leq d$, since we delete the interiors of the facets. 
Now if $\del$ is not the $n-l$-skeleton, we have $d \leq n-l$. 
Thus its codimension one
skeleton, which is $l+1$-Cohen-Macaulay, will not have maximal girth.
\remfin

The following characterizes those simplicial complexes attaining the 
submaximal $d$. Actually the characterization is maybe more transparently given
in terms of the Alexander dual $\del^*$.
Recall that a {\it missing face} $F$ is a subset of $[n]$ not in $\del$.
Also the frame dimension of $\del^*$ is denoted $c^*-1$.

\begin{proposition} \label{TopProSubm}
The following are equivalent for a simplicial complex 
$\del$.

\begin{itemize}
\item[\it a.] $\del$ is $l$-CM with $d$ equal to $n-l$.

\item[\it b.] $d$ is $n-l$ and the cardinality of $F \cup G$ is 
$\geq n+2-l$ for any two distinct minimal missing faces.

\item[\it c.] $c^*$ is $l-1$ and any two distinct facets of $\del^*$
intersect in a subset of cardinality less than $l-1$.
\end{itemize}

The girth of $\del$ being maximal $n+2-l$ corresponds to the cardinality of
any minimal missing face being $\leq n-l$, respectively 
every facet of $\del^*$ having cardinality $\geq l$.
\end{proposition}

\begin{proof} The equivalence of {\it b}. and {\it c.} 
is clear since $F$ is a minimal
missing face of $\del$ iff the complement $F^c$ is a facet of $\del^*$,
and $c^* + d + 1$ is $n$.

{\it a.} $\impl$ {\it b.} 
Assume $\del$ is $l$-CM with $d$ equal to $n-l$. Let $R$ be a subset of 
cardinality
$l-1$. Then $\del^\prime = \del_{-R}$ is CM with $d$ equal to  
$n-l =: n^\prime -1$. 
We must show that $\del^\prime$ has only one minimal missing face. But this
is true, the only minimal missing face is the intersection of the missing faces
of cardinality $n^\prime -1$.

\medskip Before proving the converse, note that if $\del$ is a simplicial 
complex containing only one minimal missing face $F$, clearly $F$ is non-empty.
Letting $x$ be in $F$, the restriction $\del_{-\{x\}}$ contains no minimal
missing face and so is a simplex and therefore $d$ is $n-1$.

{\it b.} $\impl$ {\it a.} 
Let $d$ be $n-l$ with $l \geq 1$ and $R$ have cardinality
$l-1$. Then $\del^\prime = \del_{-R}$ contains at most one minimal missing
face. Since $d^\prime \leq n-l$, by the above it contains exactly one minimal 
missing face $F$ and the facets are exactly the $n^\prime -1$-set of 
$[n]\backslash R$
not containing $F$. Therefore $\del^\prime$ is CM of the same dimension as
$\del$, and so $\del$ is $l$-CM with $d$ equal to $n-l$.  

\medskip
To prove the last statement, note that any minimal missing face of $\del$
has cardinality $\leq d+1$ which is  $n+1-l$. 
That the girth of $\del$ is maximal 
$n+2-l$ means that $\tH_{\dd}(\del_{-R})$ which is $\tH_{n-l-1}(\del_{-R})$
vanishes when $R$ is of cardinality  $l-1$. 
But then $[n]\backslash R$ is not a minimal missing
face and so all these have cardinality $\leq n-l$.
\end{proof}

\rem \label{TopRemFac}
When $l$ is $2$ we see that for the Alexander dual of a $2$-CM 
simplicial complex $\del$
with $d$ submaximal equal to  $n-2$, 
the facets partition $[n]$ into disjoint subsets.
In general for $l$-CM $\del$ with $d$ equal to $n-l$, 
the facets of its Alexander dual are
a collection of subsets $F_1, \ldots, F_m$ such that each $l-1$-subset is
contained in exactly one subset $F_i$. The girth of $\del$ is $n+2-l$
if all the $F_i$ have cardinality $\geq l$.
\remfin

If $\dl$ is $3$-CM and if $d$ is $2$, i.e. $\dl$ is a graph which is
$3$-connected, it is rather clear that if $\dl$ has a reasonably large 
number of vertices, then $\dl$ cannot have maximal girth $n-1$.
This suggests that for $l$-CM $\dl$ of maximal girth and $n$ reasonably
large, the dimension will not be to small. This would be a consequence
of the following.

\begin{conjecture} \label{TopConRam}
Let $\dl$ be an $l$-CM simplicial complex of maximal girth. Assume
it is not the $r$-skeleton of the $l+r-1$-simplex for some $r$. Then 
$c \geq l-1$.
\end{conjecture}

This of course implies $d \geq l$ if $\dl$ is not the skeleton of some
simplex of dimension $\leq 2l-3$.

We prove this conjecture in the case of $d$ submaximal equal to $n-l$. In
fact we prove something stronger.

\begin{proposition}
Let $\dl$ be an $l$-CM simplicial complex of maximal girth with $d$ equal to 
$n-l$. 
Then $c \geq l-1$ and for $l \geq 3$ and $n \geq 3l-4$ we have the stronger
bound $c \geq (n+2-l)/2$.
\end{proposition}

\begin{proof} Clearly the statement holds if $l$ is $1$ or $2$.
So we assume $l \geq 3$ and 
consider the Alexander dual simplicial complex $\dl^*$.
Let $X$ be a facet of $\dl^*$ of cardinality $d^*$, which is
$\geq l$ by Proposition \ref{TopProSubm}, and let $Y$
be  the complement $[n]\backslash X$. We shall count the number of 
pairs $(A,B)$ where :
\begin{itemize}
\item[i.] $A \sus X$ has cardinality $l-2$,
 \item[ii.] $B \sus Y$ has cardinality $2$,
\item[iii.] $A \cup B$ is a face of $\dl^*$.
\end{itemize}

First given $A$ fulfilling i. Since $c^* = l-1$, we may take any $y$ in $Y$
and $A \cup \{y \}$ will be a face. By Proposition \ref{TopProSubm}, 
it may be extended to a face $A \cup \{y,z\}$. Now $z$ cannot be in $X$
since any two facets of $\dl^*$ intersect in cardinality $\leq l-2$.
Hence $z$ is in $Y$. This gives that the number of pairs $(A,B)$ fulfilling
i., ii., and iii. is greater or equal to 
\begin{equation} \label{TopForAvink}
\frac{n-d^*}{2} \cdot \binom{d^*}{l-2}.
\end{equation}

Now given $B$ fulfilling ii. Consider the set of $A$'s fulfilling i.
and iii. If $A_1$ and $A_2$ are two distinct such, then the union of 
$A_1,A_2$ and $B$ is not a face, since it would intersect $X$ in 
cardinality $\geq l-1$. Hence $A_1 \cup B$ and $A_2 \cup B$
are in distinct facets and so $A_1$ and $A_2$ intersect in 
cardinality $\leq l-4$. So for a given $B$ the number of possible
$A$'s is less or equal to 
\begin{equation} 
 \frac{1}{l-2} \cdot \binom{d^*}{l-3}
\end{equation}
Summing over the $B$'s,
the number of pairs $(A,B)$ fulfilling i., ii., and iii.
is less or equal to 
\begin{equation} \label{TopForBvink}
 \binom{n-d^*}{2} \cdot \frac{1}{l-2} \cdot \binom{d^*}{l-3}
\end{equation}
This gives (\ref{TopForAvink}) less or equal to (\ref{TopForBvink})
which gives
\begin{equation*}
d^* -l+3   \leq  n-d^* -1. 
\end{equation*}
Putting $d^* = n-c-1$ this becomes
\[ c \geq \frac{n+2-l}{2} \]

Now let $\overline{A}$ be the complement of $A$ in $X$, of cardinality
$d^*-l+2$. For a given $B$ any 
two distinct $\overline{A_1}$ and $\overline{A_2}$
intersect in cardinality $\leq d^* -l$. Thus the number of pairs
$(\overline{A},B)$ where $(A,B)$ fulfills i., ii., and iii. is less 
or equal to 
\begin{equation} \label{TopForBCvink}
\binom{n-d^*}{2} \cdot \frac{1}{d^*-l+2} \cdot \binom{d^*}{d^*-l+1}.
\end{equation}
This gives (\ref{TopForAvink}) less or equal to  (\ref{TopForBCvink}) 
which gives
\[  1 \leq \frac{1}{d^*-l+2} \cdot \frac{d^*-l+2}{l-1} \cdot (n-d^*-1). \]
Putting $d^* = n-c-1$ this becomes
\[ c \geq l-1. \]
\end{proof}

\section{Gorenstein* complexes and orientable homology manifolds}

Another natural question concerning the top cohomology module of CM-complexes
is when it can be identified as an ideal in $S$, or more intrinsically
as a rank one torsion free module over $S$. The answer is that this happens
exactly when $\del$ is Gorenstein*. More generally we show that for
Buchsbaum complexes it happens when it is an orientable homology manifold.

\subsection{Gorenstein* complexes}
Recall that $\del$ is a {\it Gorenstein} simplicial complex iff the 
Stanley-Reisner ring $k[\del]$ is a Gorenstein ring. If $\del$ is also
not a cone, it is called {\it Gorenstein*}. By a criterion of Hochster
this latter is equivalent to $\del$ being CM and 
$\tH_p(\lkd R) = k$ whenever  $p+r = \dim \del$ and $R$ is a face of $\del$. 

\begin{theorem} \label{TopTheGor}
Let $\del$ be CM. The top cohomology module of $\del$
is a rank one torsion free $S$-module iff $\del$ is a Gorenstein* simplicial 
complex. The top cohomology module then identifies as the ideal in $S$ of
the Stanley-Reisner ring of the Alexander dual simplicial complex 
of $\del$.
\end{theorem}

\begin{proof}
Note first that if $M$ is square free then $M^\vee = \Hom_S(M,\omega_S)$
is also square free as can be seen by taking a free presentation of $M$.

Consider then a non-zero map $S(-\aa) \pil \EH^{\dim \del}(\del)^\vee$
where $S(-\aa)$ is square free. Dualized this gives
a nonzero composition
\begin{equation}
\EH^{\dim \del}(\del) \pil \EH^{\dim \del}(\del)^{\vee \vee} \pil 
S(-\aa^c) \pil S \label{GorLabKom}
\end{equation}

If the top cohomology module is torsion free of rank one this map is
an inclusion. If we take the Alexander dual we get a surjection
\begin{equation} S \pil \Ext_S^{n-d}(k[\del], \omega_S). \label{GorLabAld}
\end{equation}
The module on the right of (\ref{GorLabAld}) is a $k[\del]$-module and must 
therefore be
a quotient of $k[\del]$. 
Since $k[\del]$ is unmixed, 
the right side of (\ref{GorLabAld}) has the same asscoiated prime ideals, 
see \cite{Vo} Proposition 3.6.b on p.51.
Therefore since $k[\del]$ is reduced 
they must be equal.
Thus $k[\del]$ becomes Gorenstein and Gorenstein* since $\del$ is not acyclic.

Conversely if $\del$ is Gorenstein* the right hand of (\ref{GorLabAld})
is $k[\del]$. The Alexander dual of (\ref{GorLabAld}) which is 
the composition in (\ref{GorLabKom}) then identifies the top cohomology 
module as $I_{\del^*}$.
\end{proof}

\rem When $\del$ is zero-dimensional, i.e. a set of vertices,
$\del$ is $l$-CM if $\del$ consists of $l$ or more vertices.
Thus if $\del$ is an $l$-CM simplicial complex any link $\lkd Q$ which
is a non-empty point set will consist of $l$ or more points. When $l$ is equal
to $2$,
Gorenstein* complexes are the $2$-CM complexes where each such link
consists of exactly two points. Hence the class of $l$-CM complexes
such that the links which are non-empty point sets consist of exactly
$l$ points might be a reasonable generalization of Gorenstein*
simplicial complexes.

\begin{theorem} Let $\del$ be a Gorenstein* simplicial complex. The
homology modules $\EH_i(\del)$ and $\EH_{\dd-1-i}(\del)$ are Alexander dual
square free modules.
\end{theorem}

\begin{proof}
In the resolution of $k[\del]$, in each linear strand 
the differential is given by maps, \cite{Ho},
\begin{equation}
 S \te_k \tH^p(\del_{R\cup \{i\}}) \pil  S \te_k \tH^p(\del_{R}) 
\label{TopForGor} \end{equation}
When $\del$ is Gorenstein*, the resolution of $k[\del]$ is self-dual
and $\Hom_S(-,\omega_S)$ of (\ref{TopForGor}) identifies as
\[ S \te_k \tH^{\dim \del -1-p}(\del_{-R}) \pil S \te_k 
\tH^{\dim \del -1-p}(\del_{-R - \{i \}}) \]
Hence the natural maps 
\[\tH_p(\del_R) \pil \tH_p(\del_{R \cup \{i\}}) \]
and 
\[ \tH_{\dim \del -1-p}(\del_{-R-\{i\}}) \pil 
\tH_{\dim \del -1-p}(\del_{-R}) \]
are dual to each other. But this means that the square free modules
$\EH_i(\del)$ and $\EH_{\dim \del -1-i}(\del)$ are Alexander duals.
\end{proof}

\subsection{Orientable homology manifolds}

Recall that a connected $\del$ is an {\it orientable homology manifold} 
if all proper links
of $\del$ are Gorenstein* and $\tH_{\dd}(\del)$ is $k$. Then in particular
$\del$ is Buchsbaum.

\begin{theorem} \label{TopTheBb} 
Let $\del$ be Buchsbaum of dimension $\geq 1$. 
Then $\del$ is a connected orientable 
homology manifold iff the top cohomology module is a rank one torsion free $S$-module.
It may be identified as the ideal in $S$ of the Stanley-Reisner ring
of the Alexander dual simplicial complex of $\del$.
\end{theorem}

\begin{proof}
If the top cohomology module is torsion free of rank one we
get as in Theorem \ref{TopTheGor} an inclusion
\begin{equation}
\EH^{\dim \del}(\del) \pil S \label{BBLabKom}
\end{equation}
Taking the Alexander dual we get a surjection
\begin{equation} S \pil \Ext_S^{n-d}(k[\del], \omega_S). \label{BBLabAld}
\end{equation}
The module on the right of (\ref{BBLabAld}) is a $k[\del]$-module and as such
a quotient of $k[\del]$. Also it has the same associated prime ideals
as $k[\del]$ since the latter is unmixed, see \cite{Vo}, Proposition 3.6.b on
p.51 and Corollary 2.4 on p.75.
Since $k[\del]$ is reduced the right hand side of (\ref{BBLabAld})
must then be equal to  $k[\del]$. 
Taking the Alexander dual of (\ref{BBLabAld}) we get that the left side
of (\ref{BBLabKom}) identifies as $I_{\del^*}$.

\medskip
Since  $\EH^p(\lk_{\del}\{ x \})$ is $\EH^p(\del)_{|[n] \backslash \{x\}}$
by Proposition \ref{HomProCoh}, 
it follows by Theorem \ref{TopTheGor} that $\lk_\del \{x \}$
is Gorenstein*. Therefore $\del$ is an orientable homology manifold.

Conversely suppose $\del$ is a connected orientable homology manifold. Then
$\EH^{\dim \del}(\del)$ has rank one. If it had torsion there
would have to be a proper subset $R$ of  $[n]$ such that 
$\EH^{\dim \del}(\del)_{|R}$ 
had rank $\geq 2$. But since this module identifies as 
$\EH^{{\dim \del}-r^c} (\lk_\del R^c)$ and this link is Gorenstein*, this is
impossible. Thus the top cohomology module is torsion free of rank one.
\end{proof}

\section{Designs and vanishing of homology modules }

Assume $\del$ is not $\emptyset$ or the simplex on $[n]$. 
Recall that $c$ is the maximum integer $i$ such that all $i$-sets are in
$\del$. If $T$ is a $c+1$-set not in $\del$, then $\tH_{c-1}(\del_T)$ is
non-zero. Hence the homology module $\EH_{c-1}(\del)$ is non-zero.
The following, from \cite{FV}, describes when the other 
homology modules vanish.

\subsection{Bi-Cohen-Macaulay complexes}

\begin{theorem} \label{VanTheBi}
There is at most one non-vanishing homology module 
$\EH_i(\del)$ (with the exception of $\EH_{-1}(\del)$ if this is $k$) iff
the Alexander dual of $\del$ is Cohen-Macaulay.
\end{theorem}

\rem $\EH_{-1}(\del)$ is $k$ iff all vertices of $[n]$ are in $\del$.
\remfin

When both $\del$ and $\del^*$ are Cohen-Macaulay we call $\del$ 
{\it bi-Cohen-Macaulay.} This corresponds to $\del$ having one non-vanishing
cohomology module and one non-vanishing homology module (save the exception). 

\eks When $\del$ is a graph, $\del$ is Cohen-Macaulay iff it is connected.
The Alexander dual of $\del$ is Cohen-Macaulay iff $\del$ is a forest.
Hence $\del$ is bi-Cohen-Macaulay iff it is a tree.
\eksfin

In \cite{FV} it was shown that the $f$-vector of a bi-Cohen-Macaulay 
simplicial complex only depends on the number of vertices $n$, its
dimension, and its frame dimension. If $f_{\del}(t) = \sum f_{i-1}t^i$
is the $f$-polynomial then
\begin{equation} \label{VanForBic} 
f_{\del}(t) = (1+t)^{d-c}(1 + (n-d+c)t + \cdots +\binom{n-d+c}{c} t^c). 
\end{equation}

\medskip

\subsection{$l$-Cohen-Macaulay designs}

Our objective is now to put strong vanishing conditions on the homology
modules of $l$-CM simplicial complexes and investigate what kind of
simplicial complexes we obtain this way.

We define a simplicial complex to be an {\it $l$-CM design} iff
i) $\del_{-R}$ is bi-CM of the same dimension and frame dimension as $\del$
for all $R$ of cardinality $l-1$ and ii) $\del$ is not the $(n-p)$-skeleton
of the simplex for some $p > l$.

\begin{theorem} \label{VanTheDes}
$\del$ is an $l$-CM 
design iff $\del$ is $l$-CM, of maximal girth $n+2-l$, and has at most one
non-vanishing homology module $\EH_i(\del)$ for $i < \dim \del$ (with
the exception of $\EH_{-1}(\del)$ if this is $k$).
\end{theorem}

\begin{proof}
First 
assume $\del$ is an $l$-CM design. If $d > c$, then $\tH_{\dd}(\del_{-R})$
vanishes for $R$ of cardinality $l-1$. 
So $\EH_{\dd}(\del)$ is zero in degrees $\leq n+1-l$
and the girth of $\del$ is $n+2-l$. This is also true if $d$ 
(and hence $c$) is $n+1-l$.

Since $\del$ is $l$-CM, $\tH_i(\del_{-R})$ vanishes when $i+r \leq d+l-3$ and
$i \leq d-2$. In particular when $r \leq l-1$ and $i \leq d-2$, this vanishes. 
When $i$ is not equal to $c-1$ and $r \geq l$ it will also vanish as we now
explain.
Let $T$ be a subset of $R$ of cardinality $l-1$. $\del_{-T}$ is
bi-CM of frame dimension $c-1$, and so $\EH_i(\del_{-T})$ vanishes for 
$i$ not equal to $c-1$ (except for $i = -1$ when this homology module is
$k$). Hence $\EH_i(\del)$ vanishes for  
$i < \dim \del$, except  
for  $i$ equal to $c-1$ (and for $i= -1$ when this homology module is $k$).

\medskip
Conversely, suppose $\del$ is $l$-CM, has maximal girth $n+2-l$ and at most
one non-vanishing $\EH_i(\del)$ for $i < \dim \del$ (save the exception). 
If $\del$ is a skeleton of the simplex, the condition of maximal girth clearly
implies that $\del$ is the $n-l$-skeleton.

Let $R$ be a subset of $[n]$ of cardinality $l-1$. Then $\EH_i(\del)_\bb$
is equal to $\EH_i(\del_{-R})_\bb$ when $\bb$ has support in $[n]\backslash R$.
Thus the girth of $\del$ being maximal implies the vanishing of 
$\EH_\dd(\del_{-R})$. Also the vanishing of $\EH_i(\del)$ for $i < \dd$
except when $i$ is $c-1$ (and $i=-1$ if this homology module is $k$), 
implies the same for $\del_{-R }$. 
Hence $\del_{-R}$ is bi-CM of the same dimension and frame
dimension as $\del$.
\end{proof}



\rem A slight nuisance in the formulations of Theorems \ref{VanTheBi}
and \ref{VanTheDes} is the exception statement. It may be avoided
with the following approach.
One can sheafify the complexes $L(\del)$ to get a complex
of coherent sheaves on the projective space ${\bf P}^n$
\[ \tilde{L}(\del) : \,\, \cdots 
\pil  \gO_{{\bf P}^n}(-p) \te_k
(C_\del)_p \pil \cdots \]
and consider the homology sheaves $H_{p+1}(\tilde{L}(\del))$ instead of the
enriched homology modules $\EH_{p}(\del)$. 
Then $H_{0}(\tilde{L}(\del))$
is zero when $\EH_{-1}(\del)$ is $k$ or $0$, and non-zero otherwise, so
one avoids the exception in the theorems.
\remfin

\eks \llabel{VanEksCyc} The boundaries of cyclic polytopes of even dimension 
$d$ (and containing all vertices
of $[n]$) are examples of $2$-CM designs. Since the boundaries of cyclic 
polytopes are 
Gorenstein* they are evidently $2$-CM. Furthermore since they are 
$d/2$-neighbourly, i.e. $d$ is $2c$, $\EH_i(\del)$ vanishes for $i$ in the
interval from $0$ to $c-2$. By Alexander duality $\EH_{d-2-i}(\del)$ also
vanishes when $d-2-i$ is in the interval from $c$ to $d-2$.
Hence it is a $2$-CM design.
\eksfin

We now turn to examine the rationale behind the terminology design.
Recall that a $t-(n,k,\la)$ design is a collection of (distinct)
$k$-subsets of $[n]$, called blocks, such that each $t$-subset is
contained in exactly $\la$ blocks. 

Considering the blocks as facets of a simplicial complex $\del$ this means
that $\del$ is pure of dimension $k-1$ and $\lk_\del T$ has exactly
$\la$ facets for each $t$-subset $T$. Our $l$-CM designs will be
block designs in this meaning. But they have in fact much stronger
regularity properties as we show in the following statements.

\begin{lemma}  \label{VanLemInt} 
Let $\del$ be a simplicial complex on $[n]$. 

a. Suppose that for $x$ in $[n]$ the restrictions
$\del_{-\{x\}}$all have the same $f$-polynomial $f(t)$.
Then the $f$-polynomial $f^!$ of $\del$ is given by 
\[ f^!(t) = \sum_{i=0}^{\dim \del + 1} \frac{n}{n-i} \cdot f_{i-1} t^i. \]

b. Suppose all links $\lk_\del \{x\}$ have the same $f$-polynomial $f(t)$.
Then the $f$-polynomial $f^{\#}$ of $\del$ is given by 
\[ f^{\#}(t) = 1 + n \cdot \sum_{i=0}^{\dim \del + 1} \frac{f_{i-1}}{i+1} t^{i+1}. \]
In particular $f(t)$ is the derivative of $f^{\#}(t)$ divided by the number
of vertices.
\end{lemma}

\begin{proof}
a. We count the number of pairs $(F,x)$ where $F$ is a face of cardinality
$i$ and $x$ is not in $F$. By restricting to each $\del_{-\{x\}}$ this
is $f_{i-1} \cdot n$. By counting first the $F$'s in $\del$, these pairs
can be counted as $f^!_{i-1} \cdot (n-i)$.

b. We count the pairs $(F,x)$ where $F$ is a face of cardinality $i+1$ and
$x$ is in $F$. Considering each link $\lk_\del \{x\}$ this is 
$f_{i-1} \cdot n$.
By counting first the $F$'s, this can be counted as $f^{\#}_i \cdot (i+1)$. 
\end{proof}

\begin{corollary} \label{VanCorLink} 
Let $\del$ be an $l$-CM design and $Q \sus [n]$ with 
$q \leq l-1$. Then the $f$-vector of $\lk_\del Q$ only depends on $n,d,c,l$, and $q$.
In particular, when $Q = \emptyset$, the $f$-vector of $\del$ depends
only on $n,d,c$, and $l$.
\end{corollary}

\begin{proof}
When $Q = \emptyset$ this follows by repeatedly using 
Lemma \ref{VanLemInt} a. since all restrictions $\del_{-R}$ where $r = l-1$
have the same $f$-vector. If $Q$ is nonempty let $Q = Q^\prime \cup \{x\}$.
There is a natural exact sequence
\[ 0 \pil C_{\lk_{\del_{-\{x\}}} Q^\prime} \pil C_{\lk_\del Q^\prime} 
\pil C_{\lk_\del Q}(-1) \pil 0. \] 
The statement follows by induction, on taking Hilbert series of these modules.
\end{proof}

\begin{corollary} \label{VanCorLam}
An $l$-CM design is an $l\!-\!1-(n,d,\lambda)$ block design where
$\lambda$ is 
\begin{equation} \label{VanLabFcgl}
\binom{n-d + c+1-l}{c+1-l} \cdot \binom{c}{l-1}^{-1} \cdot \binom{d}{l-1},
\text{ for } c \geq l-1, \end{equation}
and
\begin{equation} \label{VanLabFcll}
\binom{n-d}{l-1-c}^{-1} \cdot \binom{l-1}{c} \cdot \binom{d}{l-1}, 
\text{ for } c \leq l-1.\end{equation}
\end{corollary}

\begin{proof}
a. The restriction $\del_{-R}$ when $r= l-1$ is bi-CM with invariants 
$n+1-l,d$, and
$c$ so the number of facets is by (\ref{VanForBic}) given by 
$\binom{n-d+c+1-l}{c}$ which we denote by  
$\la^{\prime \prime}$. 
When $c \geq l-1$ this may be written as
\begin{equation} \label{VanLabcgl}
\binom{n-d + c+1-l}{c+1-l} \cdot \binom{c}{l-1}^{-1} \cdot \frac{1}{(l-1)!} 
\cdot \Pi_{i=2}^l (n-d+i-l), \end{equation}
and when $c \leq l-1$ this may be written as
\begin{equation} \label{VanLabcll} 
\binom{n-d}{l-1-c}^{-1} \cdot \binom{l-1}{c} \cdot \frac{1}{(l-1)!}
\cdot \Pi_{i=2}^l (n-d+i-l). \end{equation}
By Lemma \ref{VanLemInt} a.
the number of facets of $\del_{-R}$ when $r = l-2$ is 
given by $\la^{\prime \prime} \cdot \frac{n+2-l}{n-d+2-l}$ and in general
for $r \leq l-1$ by
\[\la^{\prime\prime} \cdot \Pi_{i=2}^{l-r} \frac{n+i-l}{n-d+i-l}. \]
The expression for this when $r=0$ will be the number of facets of $\del$
and we denote it by $\lambda^\prime$.
Now by Lemma \ref{VanLemInt} b. the number of facets of $\lk_\del Q$ when
$q=1$ is $\lambda^\prime \cdot \frac{d}{n}$ and in general for 
$q \leq l-1$ it will be
\[ \lambda^\prime \cdot \Pi_{i=0}^{q-1} \frac{d-i}{n-i}. \]
When $q = l-1$ this is 
\[ \lambda = \lambda^\prime \cdot \Pi_{i=2}^l \frac{d+i-l}{n+i-l}. \]
By working out the result of these steps starting from
(\ref{VanLabcgl}) and (\ref{VanLabcll}) we obtain the statement.
\end{proof}

\eks According to Conjecture \ref{TopConRam} the minimal interesting
value of $c$ for $l$-CM designs is $l-1$. And then the minimal interesting
value for $d$ is $l$. By the above we see that $\la = l$
so we get $l\!-\!1-(n,l,l)$ designs. Such designs have been constructed
when $l=3$, \cite{MaRo}, and $l=4$, \cite{Kr}, in many cases. It
is not automatic that they are $l$-CM designs but we suspect that many
of the examples constructed are.
\eksfin

A consequence of the Corollary \ref{VanCorLam} is that it supports
Conjecture \ref{TopConRam}. 

\begin{corollary} There is a function $\sigma(l,d)$ such that if 
$\del$ is an $l$-CM design with $n \geq \sigma(l,d)$ then $c \geq l-1$.
\end{corollary}

\begin{proof}
The expression of (\ref{VanLabFcll}) is for $c \leq l-2$ equal to 
\[ \frac{d(d-1) \cdots (d+2-l)} {c! (n-d) \cdots (n-d+2+c-l)}. \]
If this is an integer then
\[ \frac{d(d-1) \cdots (d+2-l)}{n-d} \]
is an integer and so the statement follows.
\end{proof}

\subsection{(a,b)-designs} Now we introduce a class of designs
which extends the class of $l$-CM designs. 
Call a simplicial complex an $(a,b)$-design if 
\begin{equation} \label{VanLabLBA} (\lk_\del B)_{-A}
\end{equation}
is bi-CM of dimension and frame dimension equal to those of $\del$ reduced by
$b$, whenever $A$ and $B$ are disjoint subsets of $[n]$ of cardinalities
$a$ and $b$. Observe that the Alexander dual of an $(a,b)$-design is
a $(b,a)$-design, since the Alexander dual of (\ref{VanLabLBA}) is
$(\lk_{\del^*} A)_{-B}$. Note also that $l$-CM designs are 
$(l\!-\!1,0)$-designs
and so $(0,l\!-\!1)$ designs are Alexander duals of $l$-CM designs.

\eks  \llabel{VanEksEE}
The standard triangulation of the real projective plane with six
vertices,  and the triangulation of the two-dimensional torus
with seven vertices are examples of $(1,1)$-designs $\del$ with $d=3$
and $c=2$. This is because each link $\lk_\del \{x\}$ is a polygon with
$n-1$ vertices and invariants
$n-1, d-1$, and $c-1$. Thus it is an $(1,0)$-design (or $2$-CM design) with 
the right invariants. Note that
both reduced  homology groups vanish for the real projective plane 
(char $k \neq 2$) while the torus has 
$\tH_{2}(\del)$ one-dimensional and $\tH_1(\del)$ 
two-dimensional. 

It is known for which genera of orientable and non-orientable surfaces
there exists $2$-neighbourly triangulations, see \cite{Ri}, \cite{JR}.
These will be examples of $(1,1)$-designs.

In general $(1,1)$-designs with $d=2c-1$ will be homology manifolds since
each $\lk_\del \{ x\} $ gives a $(1,0)$-design (or $2$-CM design) 
with invariants $d-1$ and $c-1$ where $d-1 = 2(c-1)$. Conferring Example 
\ref{VanEksCyc} one may work out, via Euler characteristics,
that the rank of the top cohomology
module of each singlepoint link is one, and hence each of these links
are Gorenstein* by Theorems \ref{TopTheSyz} and \ref{TopTheGor}.
Such homology manifolds are related to certain extremal behaviour of
the Euler characteristic, see \cite{Ku}, \cite{No}, or \cite{Lu}
for a survey. 

\eksfin

We now show that $(a,b)$-designs have very strong regularity
properties. Not only are the number of facets of links determined,
but the complete $f$-vector of all combinations of links and restrictions
up to a certain level is determined.

\begin{lemma} Let $\del$ be an $(a,b)$-design and $R$ and $Q$ be 
disjoint subsets of $[n]$ where $r+q \leq a+b$. Then the $f$-vector
of 
\begin{equation}
(\lk_\del Q)_{-R}  \label{VanLabLQR}
\end{equation}
depends only on the numerical invariants $n,d,c,a,b,r$, and $q$.
\end{lemma}

\begin{proof}
If $Q$ and $R$ are empty this follows by repeated use of
Lemma \ref{VanLemInt}. 
Suppose then that $r \leq a$ and $q \leq b$. Then 
(\ref{VanLabLQR}) is an $(a-r,b-q)$-design
and so this follows by the case just treated. If $q > b$ then let $B \sus Q$
be a subset of cardinality $b$. Then $(\lk_\del B)_{-R}$ is a $(a-r,0)$-design.
Since (\ref{VanLabLQR}) is a further link of this we get the statement
from Corollary \ref{VanCorLink}. If $r > a$ we may reduce to the case
just considered by taking the Alexander dual.
\end{proof}

\rem Considering the case when $Q$ has cardinality $a+b$ we see that 
$(a,b)$-designs are $a+b-(n,d,\la)$ block designs for
some $\la$. By taking the link of a set of cardinality $b$ we get
an $(a,0)$-design which is an $a-(n-b,d-b,\la)$ design and with invariants
$n-b, d-b$ and $c-b$. Hence $\la$ may be determined by Corollary 
\ref{VanCorLam}.
\remfin

\eks A consequence of the above remark is that $(0,l-1)$-designs have
$\la$ equal to $\binom{n+1+c-d-l}{c+1-l}$. Hence when $c = l-1$
we get $c-(n,d,1)$ designs and these are exactly the Steiner
systems $S(c,d,n)$. They are Alexander duals of the $l$-CM designs
where $d$ is submaximal equal to $n-l$.
\eksfin

\section{Problems and conjectures}

We pose the following problems.

\begin{problem}
What are the possible $f$-vectors (or $h$-vectors) of $l$-Cohen-Macaulay
simplicial complexes?
\end{problem}

The case when $l$ is equal to $1$ is classical, see \cite[II.3.3]{St}. 
When $l \geq 2$
this is likely to be a difficult problem since any answer also would
include an answer to what the $h$-vectors of Gorenstein* simplicial complexes
are. However, any conjecture about this would be highly interesting since
it would contain as a subconjecture what the $h$-vectors of Gorenstein*
simplicial complexes are. Some investigations into this problem is 
contained in \cite{Sw}.

This problem might be more tractable if $d$ is
an extremal value.

\begin{subproblem}
What are the possible $f$-vectors (or $h$-vectors) of $l$-CM simplicial 
complexes with $d$ equal to $n-l$?
\end{subproblem}

\begin{problem}
Construct $(a,b)$-designs 
for various parameters of $n,d,c,a$, and $b$.
\end{problem}

As has been pointed out this has been done in a number of particular
cases. When $a$ and $b$ are zero 
we have the bi-Cohen-Macaulay simplicial complexes
constructed in \cite{FV}. When $a$ is $1$, $b$ is $0$ and $d$ is $2c$ 
we have the cyclic 
polytopes.
When $a$ and $b$ are $1$, triangulations of surfaces has been constructed. 
For $a=0$, low values of $b$ and $c=b$ 
many Steiner systems $S(c,d,n)$ have been
constructed.
An examination of the literature on designs will most likely reveal
numerous other cases.

\medskip
Theorem \ref{VanTheDes} and Example \ref{VanEksEE} also suggest the following.

\begin{problem} Determine the homological behaviour of $(a,b)$-designs
$\del$. For instance do the dimensions of $\tH_i(\del)$ only depend
on $n,d,c,a$, and $b$, and if so, what is it?
\end{problem}

\medskip
We also recall the following from Section \ref{TopSek}.

\medskip
\noindent {\bf Conjecture \ref{TopConRam}.}
{\it Let $\dl$ be an $l$-CM simplicial complex of maximal girth. Assume
it is not the $r$-skeleton of the $l+r-1$-simplex for some $r$. Then 
$c \geq l-1$.}

\noindent or in a weaker form

\begin{conjecture} There is an integer $\sigma(l,d)$ such
that if $\del$ is an $l-CM$ simplicial complex of maximal girth where
the number of vertices is larger than this
integer, then $c \geq l-1$.
\end{conjecture}

\end{document}